\documentclass[12pt]{article}
  \usepackage{amsfonts}
  \usepackage[square,sort&compress,comma,numbers]{natbib} 
   \def\R{\mathbb{R}}

   \def\1{{\rm I\mskip -10.5mu 1}} 
   
   \def\e{{\varepsilon}}

   \def\cC{{\cal C}}

   \def\div{\mathop{\rm div}\nolimits}

   \def\no{\noindent}
   \def\proof{\mbox {{\underline {\sf Proof}} \hspace{2mm}}}
   \def\qed{{\hfill {\em q.e.d.}\\\vspace{1mm}}}
   
   \newcommand{\beq}{\begin{equation}}
   \newcommand{\eeq}{\end{equation}}

\newtheorem{df}{Definition}[section]
\newtheorem{prop}[df]{Proposition}
\newtheorem{lemma}[df]{Lemma}
\newtheorem{teo}[df]{Theorem}

\newtheorem{cor}[df]{Corollary}

 \newcommand{\sezione}[1]{\section{#1}\setcounter{equation}{0}}

\begin{document}

%%%%%%%%%%%%%%%%%%%%%%%%        TITOLO       %%%%%%%%%%%%%%%%%%%%%%%%%%%%%

   \title{Nonexistence of solutions for elliptic equations with
     supercritical nonlinearity in nearly nontrivial domains
%\\
%\vspace{3mm}
%Titolo in francese
}

% \date{}

  \maketitle
%\vspace{-2cm}

%%%%%%%%%% quando tolgo la data devo ripristinare il vspace da -2 cm

%\vspace{-2cm}
 \vspace{5mm}

\begin{center}

{ {\bf Riccardo MOLLE$^a$,\quad Donato PASSASEO$^b$}}

\vspace{5mm}

{\em
${\phantom{1}}^a$Dipartimento di Matematica,
Universit\`a di Roma ``Tor Vergata'',\linebreak
Via della Ricerca Scientifica n. 1,
00133 Roma, Italy.}

\vspace{2mm}

{\em
${\phantom{1}}^b$Dipartimento di Matematica ``E. De Giorgi'',
  Universit\`a di Lecce,\linebreak 
P.O. Box 193, 73100 Lecce, Italy.
}
\end{center}

\vspace{5mm}

%%%%%%%%%%%%%%%%%%%%%%          ABSTRACT         %%%%%%%%%%%%%%%%%%%%%%%%%%

{\small {\sc \noindent \ \ Abstract.} - 
\footnote{ {\em E-mail address:} molle@mat.uniroma2.it (R. Molle).}
We deals with nonlinear elliptic Dirichlet problems of the form 
$$
\div(|D u|^{p-2}D u )+f(u)=0\quad\mbox{ in }\Omega,\qquad u\in
H^{1,p}_0(\Omega)
$$
where $\Omega$ is a bounded domain in $\R^n$, $n\ge 2$, $p> 1$ and
$f$ has supercritical growth from the viewpoint of Sobolev embedding.

\no Our aim is to show that there exist bounded contractible non
star-shaped domains $\Omega$, arbitrarily close to domains with
nontrivial topology, such that the problem does not have nontrivial
solutions. 
For example, we prove that if $n=2$, $1<p<2$, $f(u)=|u|^{q-2}u$ with
$q>{2p\over 2-p}$ and $\Omega=\{(\rho\cos\theta,\rho\sin\theta)\ :\
|\theta|<\alpha,\ |\rho -1|<s\}$ with $0<\alpha<\pi$ and $0<s<1$, then
for all $q>{2p\over 2-p}$ there exists $\bar s>0$ such that the
problem has only the trivial solution $u\equiv 0$ for all $\alpha\in
(0,\pi)$ and $s\in (0,\bar s)$.

\vspace{3mm}

%%%%%%%%%%%%%%%%%%%%  A.M.S. Subj. &   KEY WORDS  %%%%%%%%%%%%%%%%%

{\em  \noindent \ \ MSC:}  35J20; 35J60; 35J65.

\vspace{1mm}

{\em  \noindent \ \  Keywords:} 
   Supercritical Dirichlet problems, contractible domains,
   nonexistence of solutions.
}

%%%%%%%%%%%%%%%%%%%%%%%%     I Paragrafo     %%%%%%%%%%%%%%%%%%%%%%%%%%%%%%%%

\sezione{Introduction}

    %%%%%%%%%%%%%%%%%%%%%%%%%%%%%%%%%%%%%%%%%%%%%%%%%%%%%%%%%%%%
    %%%%%%%%                                          %%%%%%%%%%
    %%%%%%%%             I   Paragrafo                %%%%%%%%%%
    %%%%%%%%                                          %%%%%%%%%%
    %%%%%%%%%%%%%%%%%%%%%%%%%%%%%%%%%%%%%%%%%%%%%%%%%%%%%%%%%%%%

Let us consider the Dirichlet problem
\beq
\div(|D  u|^{p-2}D  u )+f(u)=0\quad\mbox{ in }\Omega,\qquad u\in
H^{1,p}_0(\Omega)
\eeq
where $\Omega$ is a bounded domain of $\R^n$, $n\ge 2$ and $p>1$.

\no It is well known that if the function $f:\R\to\R$ has critical or
supercritical growth from the viewpoint of the Sobolev embedding
$H^{1,p}_0 \hookrightarrow L^q(\Omega)$, the usual methods to find
solutions of this problem do not work (see for instance \cite{BN}).

\no For example, if $1<p<n$ and $f(t)=|t|^{q-2}t$ with $q\ge {np\over
  n-p}$ (the critical Sobolev exponent), then the existence of
nontrivial solutions to problem 
\beq
\label{*}
\div(|D  u|^{p-2}D  u )+|u|^{q-2}u=0\quad\mbox{ in }\Omega,\qquad u\in
H^{1,p}_0(\Omega)
\eeq
is strictly related to the shape of $\Omega$.
If $\Omega$ is star-shaped, problem (\ref{*}) has only the trivial
solution $u\equiv 0$, as a consequence of a Pohozaev type identity
(see \cite{Po}). 
On the other hand, if $\Omega$ is an annulus, one can easily find
infinitely many radial solutions (as pointed out by Kazdan and Werner
in \cite{KW}).
Hence, many researches have been devoted to study the effect of the
domain shape on the existence of nontrivial solutions for problem
(\ref{*}), following some stimulating questions posed by Brezis,
Nirenberg, Rabinowitz, etc. \ldots (see \cite{B}).
In particular, the case where $p=2$, $n\ge 3$, $q\ge {2n\over n-2}$
has been considered in many papers.

\no Answering a question of Nirenberg, Bahri and Coron proved in
\cite{BaCo} the existence of a positive solution when $p=2$, $n\ge 3$,
$q={2n\over n-2}$ and $\Omega$ has nontrivial topology, in the sense
that some homology group is nontrivial (see also \cite{Co,R},
concerning the case of domains with small holes). 

\no Notice that for $q>{2n\over n-2}$ the condition that $\Omega$ has
nontrivial topology is neither sufficient nor necessary to guarantee
the existence of nontrivial solution.
In fact (answering a question posed by Rabinowitz) the second author
proved  in \cite{Pjfa93,Pdie95}  that there exist exponents
$q>{2n\over n-2}$ and nontrivial domains $\Omega\subset \R^n$ with
$n\ge 3$ such that the problem  
\beq
\label{q}
\Delta u+|u|^{q-2}u=0\ \mbox{ in }\Omega,\qquad u= 0\ \mbox{ on }\partial \Omega
\eeq
has only the trivial solution $u\equiv 0$.

\no Moreover, for all $q\ge {2n\over n-2}$ there exist contractible
domains $\Omega \subset \R^n$ with $n\ge 3$ such that problem
(\ref{q}) has positive and sign-changing solutions (see
\cite{Di,D88,Pmm89,P93,P94,Pl92,P4092,Ptmna96,Pd98,MPcvpde06,MPaihp06,
MPcras02,MPcras2002} and the references therein). 

\no More precisely, for all $\alpha\in(0,\pi)$ and $s\in(0,1)$, let us
consider for example the piecewise smooth contractible domain $\Omega$
of the form 
\beq
\label{Omega}
\Omega_n^{\alpha,s}=\{(x,y)\in\R\times\R^{n-1}\ :\ (x,|y|)\in S^{\alpha,s}\}
\eeq
where
\beq
\label{s}
S^{\alpha,s}=\{(\rho\cos\theta,\rho\sin\theta)\in\R^2\ :\
0\le\theta<\alpha,\ |\rho-1|<s\}. 
\eeq
Then, the following assertions hold for problem (\ref{q}) with
$\Omega=\Omega_n^{\alpha,s}$ and $n\ge 3$: 
\begin{itemize}
\item
for all $q\ge {2n\over n-2}$ there exists $\bar \e_q>0$ such that, if
$\pi-\bar\e_q<\alpha<\pi$, then problem (\ref{q}) has positive and
sign changing solutions; moreover, for $q>{2n\over n-2}$, the number
of solutions tends to infinity as $\alpha\to \pi$ (see
\cite{Pmm89,Pd98,P4092,Ptmna96,P94,MPaihp06,MPcvpde06,MPaihp04},
etc. \ldots); 
\item for all $\alpha>{\pi\over 2}$ there exists $\bar q_\alpha\ge
  {2n\over n-2}$ such that problem (\ref{q}) with $q\ge \bar q_\alpha$
  has at least one positive solution (see \cite{MPcvpde06}); 
\item for all $\alpha>{\pi\over 2}$ there exists $\bar\e_\alpha>0$
  such that problem (\ref{q}) with ${2n\over n-2}<q<{2n\over
    n-2}+\bar\e_\alpha$ has positive solutions (see
  \cite{MPcras02,MPaihp04}, etc. \ldots). 
\end{itemize} 
These results (that have been stimulated by an interesting question
posed by Brezis in \cite{B}) show that, even if the Pohozaev
nonexistence result can be extended to non star-shaped domains (see
\cite{CCL,DZ} and also \cite{PT,MPT,Pna95} for related phenomena), it
cannot be extended to all contractible domains when $p=2$ and $n\ge
3$. 

\no The nonexistence result obtained in the present paper, on the
contrary, suggests that the situation is quite different if $n=2$ and
$1<p<2$. 
In fact, as a direct consequence of Theorem \ref{T2.4}, we have the
following proposition. 

\begin{prop}
\label{P1.1}
Assume $n=2$ and $1<p<2$. 
Then, for all $q>{2p\over 2-p}$ there exists $\bar s\in (0,1)$ such
that problem (\ref{*}) with $\Omega=\Omega_2^{\alpha,s}$ has only the trivial
solution $u\equiv 0$ for all the pairs $(\alpha,s)$ such that $s\in (0,\bar
s)$ and $\alpha\in (0,\pi)$. 
\end{prop}

\no Since, for all $s\in (0,\bar s)$, the domain $\Omega_2^{\alpha,s}$
is contractible for all $\alpha\in(0,\pi)$,  is star-shaped for
$\alpha$ small enough and is close to a domain with nontrivial
topology when $\alpha$ is close to $\pi$, Proposition \ref{P1.1}
suggests the following natural question (analogous to the well known
one posed by Brezis in \cite{B}): if $n=2$ and $1<p<2$, can one extend 
Pohozaev's nonexistence result for star-shaped domains to all the 
contractible domains of $\R^2$? 

\no The nonexistence result presented in this paper suggests that this
question might have a positive answer. 

%%%%%%%%%%%%%%%%%%%%%%%%     II Paragrafo     %%%%%%%%%%%%%%%%%%%%%%%%%%%%%%%%

\sezione{Integral identity and nonexistence result}

    %%%%%%%%%%%%%%%%%%%%%%%%%%%%%%%%%%%%%%%%%%%%%%%%%%%%%%%%%%%%
    %%%%%%%%                                          %%%%%%%%%%
    %%%%%%%%            II   Paragrafo                %%%%%%%%%%
    %%%%%%%%                                          %%%%%%%%%%
    %%%%%%%%%%%%%%%%%%%%%%%%%%%%%%%%%%%%%%%%%%%%%%%%%%%%%%%%%%%%

The following lemma generalizes Pohozaev identity.

\begin{lemma}
\label{L2.1}
Let $\Omega$ be a piecewise smooth bounded domain in $\R^n$, $n\ge 2$
and $p>1$. 
Assume that $u\in H^{1,p}_0(\Omega)$ is a solution of the equation
\beq
\label{eq}
\div(|Du|^{p-2}Du)+f(u)=0\quad\mbox{ in }\Omega,
\eeq
where $f:\R\to\R$ is a continuous function.

\no Then, for all $v=(v_1,\ldots,v_n)\in \cC^1(\overline\Omega,
\R^n)$, the function $u$ satisfies the integral identity 
\beq
\label{id}
\left(1-{1\over
    p}\right)\int_{\partial\Omega}\hspace{-1mm}|Du|^pv\cdot \nu\,
d\sigma=\int_\Omega \hspace{-1mm}|Du|^{p-2}\left(dv[Du]\cdot
  Du\right)dx+\int_\Omega\hspace{-1mm}\div v\left[F(u)-{1\over
    p}|Du|^p\right] dx, 
\eeq
where $\nu$ denotes the outward normal to $\partial\Omega$,
$dv[\xi]=\sum_{i=1}^n(D_iv)\, \xi_i$ $\forall
\xi=(\xi_1,\ldots,\xi_n)\in\R^n$ and $F(t)=\int_0^tf(\tau)d\tau$
$\forall t\in\R$. 

\end{lemma}
  
  \proof In order to prove (\ref{id}) it suffices to apply the
  Gauss-Green formula to the function $(v\cdot Du)
  |Du|^{p-2}Du$. 
  
  \no Thus, we obtain 
  $$
  \hspace{-5cm}
  \int_{\partial\Omega}(v\cdot Du)|Du|^{p-2}(Du\cdot\nu)d\sigma =
$$
  \begin{eqnarray}
\label{ea}
  & = &
  \int_\Omega\sum_{i=1}^nD_i\left[\sum_{j=1}^nv_jD_ju\cdot |Du|^{p-2}D_iu\right]dx
  \\
  & = &
\nonumber
\int_\Omega\sum_{i,j=1}^n\Big[D_iv_jD_ju|Du|^{p-2}D_iu+v_jD_{i,j}u|Du|^{p-2}D_iu\\
\nonumber
& & \hspace{20mm} +v_jD_juD_i(|Du|^{p-2}D_iu)\Big]dx.
\end{eqnarray}
Since $u\equiv 0$ on $\partial\Omega$, we have $Du=(Du\cdot\nu)\, \nu
$ and, as a consequence, 
\beq
\label{eb}
\int_{\partial\Omega}(v\cdot
Du)|Du|^{p-2}(Du\cdot\nu)d\sigma=\int_{\partial\Omega}|Du|^p(v\cdot
\nu)d\sigma. 
\eeq
Notice that
\begin{eqnarray}
\nonumber
\int_\Omega \sum_{i,j=1}^nv_jD_{i,j}u|Du|^{p-2}D_iu\,dx &=&
{1\over 2}\int_{\Omega}\sum_{i,j=1}^nv_j|Du|^{p-2}D_j|D_iu|^2dx\\
\label{ec} & = & {1\over p}\int_{\Omega} \sum_{j=1}^nv_jD_j|Du|^pdx
\\
\nonumber & = & {1\over p}\int_{\partial\Omega}|Du|^pv\cdot\nu\,
d\sigma-{1\over p}\int_\Omega \div v\, |Du|^pdx. 
\end{eqnarray}
Moreover, since $u$ solves equation (\ref{eq}),
\begin{eqnarray}
\nonumber \int_\Omega \sum_{i,j=1}^nv_jD_juD_i(|Du|^{p-2}D_iu)dx & = &
-\int_\Omega\sum_{j=1}^nv_jD_ju\, f(u) dx 
 \\
 \label{ed}
  =  -\int_\Omega \sum_{j=1}^nv_jD_jF(u) dx
& = & \int_\Omega\div v\cdot F(u).
\end{eqnarray}
Then, (\ref{id}) follows easily from (\ref{ea}), (\ref{eb}),
(\ref{ec}), (\ref{ed}). 

\qed
 
 \begin{lemma}
 \label{L2.2}
 On the piecewise smooth domain
 $\Omega_2^{\alpha,s}=\{(\rho\cos\theta,\rho\sin\theta)\in\R^2$ : $
 |\theta|<\alpha,\ |\rho-1|<s\}$ let us consider the vector field
 $v\in\cC^1(\overline\Omega_2^{\alpha,s},\R^2)$ defined by 
\beq
\label{v}
v(\rho\cos\theta,\rho\sin\theta)=
(\rho-1)(\cos\theta,\sin\theta)+\rho\theta(-\sin\theta,\cos\theta). 
\eeq
Then, 
\begin{itemize}
\item[a)] $v\cdot\nu>0$ on  $\partial\Omega_2^{\alpha,s}$
  $\forall\alpha\in(0,\pi)$, $\forall s\in(0,1)$; 
\item[b)] $\div v(\rho\cos\theta,\rho\sin\theta)=3-{1\over\rho}$\ 
  $\forall (\rho\cos\theta,\rho\sin\theta)\in\Omega_2^{\alpha,s}$; 
\item[c)]
  $dv(\rho\cos\theta,\rho\sin\theta)[\xi]\cdot\xi=\xi_N^2+\left(2-{1\over
      \rho}\right)\xi_T^2$\ $\forall\xi=(\xi_1,\xi_2)\in\R^2$,  
\end{itemize}
where 
\beq
\xi_N=\xi_1\cos\theta+\xi_2\sin\theta\quad\mbox{ and }\quad
\xi_T=-\xi_1\sin\theta+\xi_2\cos\theta. 
\eeq
\end{lemma}

\proof Property $(a)$ is a simple consequence of the definition of
$\Omega_2^{\alpha,s}$ and $v$. 

\no In order to prove $(b)$ and $(c)$ it suffices to notice that
\beq
\begin{array}{rcl}
dv(\rho\cos\theta,\rho\sin\theta)[\xi]&
=&\xi_N(\cos\theta,\sin\theta)+\xi_N\theta(-\sin\theta,\cos\theta)\\ 
& & -
\xi_T\theta(\cos\theta,\sin\theta)+
\xi_T\left(2-{1\over\rho}\right)(-\sin\theta,\cos\theta),
\end{array}
\eeq
as one can verify by direct computation.

\qed

\begin{cor}
\label{C2.3}
Let $\Omega=\Omega_2^{\alpha,s}$ and
$v\in\cC^1({\overline\Omega}_2^{\alpha,s},\R^2)$ be as in Lemma
\ref{L2.2}. 
Let $f$ and $F$ be as in Lemma \ref{L2.1}.
Then every solution of the Dirichlet problem 
\beq
\label{D}
\div (|Du|^{p-2}Du)+f(u)=0\quad\mbox{ in } \Omega_2^{\alpha,s},\quad
u\in H^{1,p}_0(\Omega_2^{\alpha,s}) 
\eeq
satisfies the inequality
\beq
\label{in}
0\le \left[1-{2\over p}+\left(1+{1\over p}\right)\,{s\over
    1-s}\right]\int_{\Omega_2^{\alpha,s}}|Du|^pdx+\int_{\Omega_2^{\alpha,s}}\div
v\cdot F(u)dx. 
\eeq
\end{cor}

\no The proof follows directly from Lemmas \ref{L2.1} and \ref{L2.2}
(taking into account that $\left|1-{1\over\rho}\right| \le {s\over
  1-s}$  $\forall\rho\in(1-s,1+s)$). 

\no Now, we can prove a nonexistence result for nontrivial solutions in
the domain $\Omega_2^{\alpha,s}$. 

\begin{teo}
\label{T2.4}
Let $\Omega=\Omega_2^{\alpha,s}$ be as in Lemma \ref{L2.2}, $f$ and
$F$ be as in Lemma \ref{L2.1} and assume that $1<p<2$ and there exists
$q>{2p\over 2-p}$ such that 
\beq
\label{eF} 
tf(t)\ge qF(t)\ge 0\qquad\forall t\in\R.
\eeq
Then, there exists $\bar s\in(0,1)$ such that the Dirichlet problem
(\ref{D}) has only the solution $u\equiv 0$ for every pair
$(\alpha,s)$ such that $s\in(0,\bar s)$ and $\alpha\in (0,\pi)$. 
\end{teo}

\proof Notice that $u\equiv 0$ is obviously a solution of Problem
(\ref{D}) because the assumption (\ref{eF}) clearly implies $f(0)=0$. 
Let us prove that it is the unique solution.

\no Since $0\le F(u)\le {1\over q}u f(u)$, from Lemma \ref{L2.2} and
Corollary \ref{C2.3} we obtain that every solution $u$ of the
Dirichlet problem (\ref{D}) must satisfy 
\beq
0\le \left[1-{2\over p}+\left(1+{1\over p}\right)\,{s\over
    1-s}\right]\int_{\Omega_2^{\alpha,s}}|Du|^pdx+ 
\left[2+{s\over 1-s}\right]\, {1\over q}\int_{\Omega_2^{\alpha,s}} u\, f(u) dx.
\eeq
Notice that
\beq
\int_{\Omega_2^{\alpha,s}}uf(u)dx=\int_{\Omega_2^{\alpha,s}}|Du|^pdx
\eeq
as $u$ solves the Dirichlet problem (\ref{D}).
Therefore we obtain,
 \beq
0\le\left[ 1-{2\over p}+{2\over q}+\left(1+{1\over p}+{1\over
      q}\right)\,{s\over 1-s} \right]\int_{\Omega_2^{\alpha,s}}
|Du|^pdx. 
\eeq 
Since $1-{2\over p}+{2\over q} <0$ for $q>{2p\over 2-p}$,
there exists $\bar s\in(0,1)$ such that   
$ 1-{2\over p}+{2\over q}+\left(1+{1\over p}+{1\over
    q}\right)\,{s\over 1-s} <0$ $\forall s\in(0,\bar s)$. 
Therefore, if $s\in(0,\bar s)$ and $u$ solves the Dirichlet problem
(\ref{D}), we must have  
\beq
\int_{\Omega_2^{\alpha,s}} |Du|^pdx=0,
\eeq
so the proof is complete.

\qed
 
\no Finally, notice that we obtain in particular Proposition
\ref{P1.1} when in Theorem \ref{T2.4} we choose $f(u)=|u|^{q-2}u$
(which obviously satisfies condition (\ref{eF})). 

\vspace{5mm}

%%%%%%%%%%%%%%%%%%%%%%%%     Acknowledgement   %%%%%%%%%%%%%%%%%%%%%%%%%%%%%%

{\small {\bf Acknowledgement}. The authors have been supported by the ``Gruppo
Nazionale per l'Analisi Matematica, la Probabilit\`a e le loro
Applicazioni (GNAMPA)'' of the {\em Istituto Nazio\-nale di Alta Matematica
(INdAM)} - Project: Equazioni di Schrodinger nonlineari: soluzioni con
indice di Morse alto o infinito. 

The second author acknowledges also the MIUR Excellence Department
Project awarded to the Department of Mathematics, University of Rome
Tor Vergata, CUP E83C18000100006 
}

%%%%%%%%%%%%%%%%%%%%%%%%    Bibliografia    %%%%%%%%%%%%%%%%%%%%%%%%%%%%%%

%%%%%%%%%%%%%%%%%%%%%%%%%%%%%%%%%%%%%%%%%%%%%%%%%%%%%%%%%%%%%%%%%%%%
%%%%%%%%%%%%%%%%%%%%%%%%%%%%%%%%%%%%%%%%%%%%%%%%%%%%%%%%%%%%%%%%%%%%

\end{document}